\documentclass{amsart}
\usepackage[T1]{fontenc}
\usepackage{lmodern}
\usepackage{amsmath,amssymb,amsthm,mathtools}
\usepackage[expansion=false]{microtype}
\usepackage{enumitem}
\usepackage{booktabs}
\usepackage{xcolor}
\usepackage{amsmath}

\usepackage[hidelinks]{hyperref}
\hypersetup{pdftitle={Global well-posedness for the energy-critical focusing NLS on T3},
 pdfauthor={Meng Chen and Han Wang},
 pdfsubject={Global well-posedness below the modified ground-state threshold},
 pdfkeywords={nonlinear Schrodinger equation, periodic, focusing, energy critical}}
\setlist[enumerate]{leftmargin=*,itemsep=3pt}
\allowdisplaybreaks[2]
\numberwithin{equation}{section}
\newtheorem{theorem}{Theorem}[section]
\newtheorem{proposition}[theorem]{Proposition}
\newtheorem{lemma}[theorem]{Lemma}
\theoremstyle{remark}
\newtheorem{remark}[theorem]{Remark}
\newcommand{\R}{\mathbb R}
\newcommand{\T}{\mathbb T}
\newcommand{\Z}{\mathbb Z}
\newcommand{\e}{\mathrm e}
\newcommand{\ii}{\mathrm i}
\newcommand{\dd}{\,\mathrm dx}
\newcommand{\F}{F}
\newcommand{\Ec}{E_c}
\newcommand{\Ac}{A_c}
\newcommand{\Kc}{K_c}
\newcommand{\ER}{E_{\R^3}}
\newcommand{\ET}{E_{\T^3}}

\title[Focusing energy-critical NLS on $\mathbb T^3$]
{Global well-posedness for the energy-critical focusing NLS on $\mathbb T^3$}
\author{Meng Chen}
\author{Han Wang}
\address[Meng Chen, Han Wang]{Department of Mathematics and Statistics, Beijing Institute of Technology\\
Key Laboratory of Algebraic Lie Theory and Analysis of Ministry of Education\\
Beijing 100081, China}
\email[Meng Chen]{mchen@bit.edu.cn}
\email[Han Wang]{wang\_han@bit.edu.cn}
\date{}

\begin{document}
\raggedbottom
\begin{abstract}
We prove global well-posedness for the focusing quintic nonlinear
Schr\"odinger equation on the standard three-dimensional torus for
initial data below the mass-corrected energy and Sobolev-norm
thresholds determined by the Euclidean ground state. 
This gives a subthreshold focusing counterpart to the
defocusing theorem of Ionescu and Pausader. The proof combines the
periodic profile analysis with the recent nonradial Euclidean scattering
theorem of Su and Zhao.
The mass correction yields an invariant region and a nonnegative
energy decomposition, placing every Euclidean profile below both
ground-state thresholds. Uniform bounds on finite intervals follow
from a contradiction argument on shrinking time intervals.
\end{abstract}
\maketitle

\section{Introduction}

Let $\T^3=(\R/2\pi\Z)^3$. In this paper, we consider the energy-critical focusing equation on the periodic space
\begin{equation}\label{eq:nls}
 \ii\partial_tu+\Delta u=-|u|^4u,
 \qquad u(0)=u_0\in H^1(\T^3).
\end{equation}
The mass and energy are
\begin{equation}\label{eq:conserved}
 M(f)=\|f\|_2^2,
 \qquad
 \ET(f)=\frac{1}{2}\|\nabla f\|_2^2-\frac{1}{6}\|f\|_6^6.
\end{equation}
We use the same energy formula on $\R^3$ and denote it by $\ER$.

Write
\begin{equation}\label{eq:ground}
 W(x)=\left(1+\frac{|x|^2}{3}\right)^{-\frac{1}{2}},
 \qquad
 G=\|\nabla W\|_{L^2(\R^3)}^2=\|W\|_{L^6(\R^3)}^6.
\end{equation}
Then $-\Delta W=W^5$, $\ER(W)=\frac{G}{3}$, and the sharp Euclidean Sobolev
inequality \cite{Aubin,Talenti} is
\begin{equation}\label{eq:euclidean-sobolev}
 \|f\|_{L^6(\R^3)}^2\le G^{-\frac{2}{3}}\|\nabla f\|_{L^2(\R^3)}^2.
\end{equation}
On the torus, the constant mode prevents a homogeneous Sobolev
inequality. Following Yue~\cite{Yue}, we retain the sharp Euclidean
constant and add a lower-order mass term. Define
\begin{equation}\label{eq:c0}
 c_0=\inf\left\{c\ge0:
 \|f\|_6^2\le G^{-\frac{2}{3}}
 \bigl(\|\nabla f\|_2^2+c\|f\|_2^2\bigr)
 \text{ for every }f\in H^1(\T^3)\right\}.
\end{equation}
Lemma~\ref{lem:sobolev} shows that $0<c_0<\infty$ and that the inequality
holds at $c=c_0$. For $c\ge c_0$, put
\begin{equation}\label{eq:modified}
 \Ac(f)=\|\nabla f\|_2^2+cM(f),
 \qquad
 \Ec(f)=\ET(f)+\frac{c}{2}M(f)
        =\frac{1}{2}\Ac(f)-\frac{1}{6}\|f\|_6^6.
\end{equation}
The quantity $\Ec$ is conserved.

We now state our main result. 
\begin{theorem}\label{thm:main}
Fix $c\ge c_0$. If $u_0\in H^1(\T^3)$
satisfies
\begin{equation}\label{eq:threshold}
 \Ac(u_0)<G,
 \qquad \Ec(u_0)<\frac{G}{3},
\end{equation}
then \eqref{eq:nls} has a unique global solution in
$C(\R;H^1(\T^3))\cap X^1_{\mathrm{loc}}(\R)$.
For each $0<T<\infty$, the mapping $u_0\mapsto u$ is continuous from the set of
initial data satisfying \eqref{eq:threshold}, equipped with the $H^1$
topology, into $X^1([-T,T])$.
Moreover, for every $0<E_0<\frac{G}{3}$ and $0<T<\infty$, there is a constant
$C(E_0,T,c)$ such that
\begin{equation}\label{eq:uniform-main}
 \|u\|_{X^1([-T,T])}\le C(E_0,T,c)
\end{equation}
whenever $\Ac(u_0)<G$ and $\Ec(u_0)\le E_0$.
Here $X^1$ is the standard atomic solution space recalled in
Section~\ref{sec:periodic}.
\end{theorem}

The theorem applies, in particular, with $c=c_0$.
The conditions in \eqref{eq:threshold} are the three-dimensional counterparts
of the modified Sobolev-norm and energy thresholds of Yue~\cite{Yue}.
Lemma~\ref{lem:nehari} identifies $\frac{G}{3}$ as the Nehari level of $\Ec$
and characterizes \eqref{eq:threshold} in terms of the associated Nehari
functional. Global existence is the appropriate conclusion
on the torus: even a nonzero spatially constant solution
$u(t,x)=a\e^{\ii|a|^4t}$ does not scatter to the linear periodic flow.

Concentration at small spatial scales links the energy-critical
periodic problem to its Euclidean counterpart.
Herr, Tataru, and Tzvetkov~\cite{HTT} established the critical local
theory on $\T^3$. Ionescu and Pausader~\cite{IP} combined this theory
with the Euclidean scattering theorem of Colliander, Keel, Staffilani,
Takaoka, and Tao~\cite{CKSTT} to prove global well-posedness for the
defocusing equation. Yue~\cite{Yue} obtained the four-dimensional
focusing theorem under modified ground-state thresholds.
Kwak~\cite{Kwak} treated $d\ge4$, including the nonalgebraic
nonlinearities for $d\ge5$, and left $d=3$ open in Remark~1.4.
Table~\ref{tab:periodic-gwp} summarizes these results on standard tori.

\begin{table}[htbp]
\centering
\setlength{\captionindent}{0pt}
\setlength{\belowcaptionskip}{6pt}
\caption{Global well-posedness for the energy-critical NLS on standard tori.}
\label{tab:periodic-gwp}
\small
\renewcommand{\arraystretch}{1.1}
\begin{tabular*}{\textwidth}{@{\extracolsep{\fill}}lll@{}}
\toprule
\shortstack[l]{Dimension\\\vphantom{(arbitrary $H^1$ data)}}
& \shortstack[l]{Defocusing\\(arbitrary $H^1$ data)}
& \shortstack[l]{Focusing\\(below both modified ground-state thresholds)} \\
\midrule
$d=3$ & Ionescu--Pausader~\cite{IP}
      & This paper (Theorem~\ref{thm:main}) \\
$d=4$ & Yue~\cite{Yue} & Yue~\cite{Yue} \\
$d\ge5$ & Kwak~\cite{Kwak} & Kwak~\cite{Kwak} \\
\bottomrule
\end{tabular*}
\end{table}

For the Euclidean focusing equation, Kenig and Merle~\cite{KM}
established the radial threshold theory in dimensions three through
five. For general initial data, Killip and Visan~\cite{KV} proved the
scattering theorem for $d\ge5$, and Dodson~\cite{Dodson} treated $d=4$.
In a recent breakthrough, Su and Zhao~\cite{SZ} proved scattering
for general $\dot H^1(\R^3)$ data below both ground-state thresholds.
This supplies the Euclidean input needed here, as anticipated in
\cite[Remark 1.3]{SZ}. See also Chung and Han~\cite{CH} for another
recent treatment of the three-dimensional Euclidean problem.

\begin{theorem}[Su--Zhao~{\cite[Theorem 1.1]{SZ}}]\label{thm:euclidean}
If $v_0\in\dot H^1(\R^3)$ satisfies $\ER(v_0)<\frac{G}{3}$ and $\|\nabla v_0\|_2^2<G,$ 
then the solution of $\ii\partial_tv+\Delta v=-|v|^4v$ is global,
belongs to $L^{10}_{t,x}(\R\times\R^3)$, and scatters in $\dot H^1$
in both time directions.
\end{theorem}

The polynomial quintic nonlinearity
allows us to use the classical profile method of \cite{IP}.
The Euclidean input requires both an energy bound and a gradient
bound for each profile. A periodic gradient bound alone does not
verify the energy condition. We obtain a nonnegative decomposition
of the modified focusing energy, including profiles prescribed by
scattering data, and thereby place every profile below both thresholds.

The argument uses only the finite scattering norm of each fixed
Euclidean solution supplied by Theorem~\ref{thm:euclidean}.
For the periodic background, local existence suffices.
We estimate the finitely many large Euclidean profiles individually
and control the remaining profiles by small-data theory and
orthogonality.
This gives a direct contradiction on shrinking time intervals,
without an additional induction on the periodic energy or a uniform
Euclidean scattering estimate.

On waveguides, Yu, Yue, and Zhao~\cite{YYZ} proved global
well-posedness below modified ground-state thresholds for the focusing
cubic equation on $\R\times\T^3$, a four-dimensional energy-critical
problem. For the quintic equation on waveguides of total dimension three,
Theorem~\ref{thm:euclidean} supplies the corresponding Euclidean
scattering input. Luo~\cite{Luo} had already identified this input as the
missing ingredient in extending his subthreshold scattering theorem to
$\R^2\times\T$. We expect it to be useful for further waveguide results,
with the variational thresholds and the remaining profile analysis
adapted to the geometry.

We establish the variational bounds in Section~\ref{sec:variational},
record the periodic profile estimates in Section~\ref{sec:periodic},
and prove Theorem~\ref{thm:main} in Section~\ref{sec:proof}.

\paragraph{Notation.}
We write
$\F(z)=-|z|^4z$ and $\langle\xi\rangle=(1+|\xi|^2)^{\frac{1}{2}}$.
The operators $P_N$ are smooth dyadic Fourier projections,
with $N\in\{1,2,4,\ldots\}$ and the zero frequency included in $P_1$. 
 For every $a>0$, we set 
\[P_{\le a}
 :=
 \sum_{\substack{N\in\{1,2,4,\ldots\}\\ N\le a}}P_N,
 \qquad
 P_{>a}
 :=
 \operatorname{Id}-P_{\le a},
 \]
 where an empty sum is understood as zero. 
A solution is understood in the Duhamel sense in
$C_tH^1\cap X^1_{\mathrm{loc}}$.

\section{Variational bounds}\label{sec:variational}

Throughout the remainder of the paper, we fix $c\ge c_0$.
We first establish the periodic Sobolev inequality with the
sharp Euclidean leading constant. We then identify the associated
variational threshold and derive the energy trapping estimate
used below.

\begin{lemma}\label{lem:sobolev}
The constant $c_0$ in \eqref{eq:c0} is positive and finite. For every
$c\ge c_0$ and $f\in H^1(\T^3)$,
\begin{equation}\label{eq:periodic-sobolev}
 \|f\|_6^6\le G^{-2}\Ac(f)^3.
\end{equation}
\end{lemma}

\begin{proof}
Choose a finite smooth partition $\sum_\nu\chi_\nu^2=1$ subordinate to
flat coordinate charts. Apply \eqref{eq:euclidean-sobolev} to each
$\chi_\nu f$, extended by zero in its Euclidean chart. The triangle
inequality in $L^3$ and the identity
$\sum_\nu\chi_\nu\nabla\chi_\nu=0$ give
\begin{align*}
 \|f\|_6^2
 &\le\sum_\nu\|\chi_\nu f\|_6^2
 \le G^{-\frac{2}{3}}\sum_\nu\|\nabla(\chi_\nu f)\|_2^2\\
 &=G^{-\frac{2}{3}}\Big(\|\nabla f\|_2^2+
 \int_{\T^3}\sum_\nu|\nabla\chi_\nu|^2|f|^2\dd\Big).
\end{align*}
This proves finiteness.
Since $|\T^3|=(2\pi)^3$,
testing the inequality  with $f=1$ gives
$c_0\ge \frac{G^{\frac{2}{3}}}{4\pi^2}>0$.
Letting admissible constants tend to $c_0$ proves the inequality at $c_0$;
cubing it gives
\eqref{eq:periodic-sobolev}.
\end{proof}

\begin{lemma}[Variational threshold]\label{lem:nehari}
Set $\Kc(f)=\Ac(f)-\|f\|_6^6$. Then
\begin{equation}\label{eq:nehari}
 \inf\{\Ec(f):f\in H^1(\T^3)\setminus\{0\},\ \Kc(f)=0\}
 =\frac{G}{3}.
\end{equation}
If $\Ec(f)<\frac{G}{3}$, then
\begin{equation}\label{eq:positive-well}
 \Ac(f)<G\quad\Longleftrightarrow\quad\Kc(f)\ge0.
\end{equation}
For nonzero $f$ satisfying these conditions, $\Kc(f)>0$.
\end{lemma}

\begin{proof}
If $f\ne0$ and $\Kc(f)=0$, then
$\Ac(f)=\|f\|_6^6\le G^{-2}\Ac(f)^3$.
Thus $\Ac(f)\ge G$ and $\Ec(f)=\frac{1}{3}\Ac(f)\ge\frac{G}{3}$.
For the reverse bound, let $w_N$ be a smooth cutoff of
$N^{\frac{1}{2}}W(Nx)$ in a fixed coordinate ball, with cutoff equal to
one near the origin. Its mass is $O(N^{-1})$, and direct scaling gives
$\Ac(w_N)\to G$ and $\|w_N\|_6^6\to G$.
Set $\lambda_N=\bigl(\frac{\Ac(w_N)}{\|w_N\|_6^6}\bigr)^{\frac{1}{4}}$.
Then $\Kc(\lambda_Nw_N)=0$ and
$\Ec(\lambda_Nw_N)=\frac{1}{3}\lambda_N^2\Ac(w_N)\to\frac{G}{3}$.
This proves \eqref{eq:nehari}.
Finally,
\[
 \Ec(f)=\frac{1}{3}\Ac(f)+\frac{1}{6}\Kc(f),
 \qquad
 \Kc(f)\ge\Ac(f)\left(1-\frac{\Ac(f)^2}{G^2}\right).
\]
The first identity proves one implication in \eqref{eq:positive-well},
and the second proves the other and the strict positivity assertion.
\end{proof}

\begin{lemma}[Energy trapping]\label{lem:trapping}
Suppose $\Ac(u_0)<G$ and $\Ec(u_0)\le E_0<\frac{G}{3}$. Then, throughout the
maximal lifespan of the solution to \eqref{eq:nls},
\begin{equation}\label{eq:trapping}
\Ac(u(t))\le3\Ec(u_0)\le3E_0<G.
\end{equation}
In particular, the solution is uniformly bounded in $H^1$ on its
lifespan.
\end{lemma}

\begin{proof}
Set $h(q)=\frac{q}{2}-\frac{q^3}{6G^2}$. By
\eqref{eq:periodic-sobolev}, $\Ec(u)\ge h(\Ac(u))$.
The function $h$ is strictly increasing on $[0,G)$ and $h(G)=\frac{G}{3}$.
Conservation of $\Ec$ and continuity of $\Ac(u(t))$ therefore prevent
$\Ac(u(t))$ from reaching $G$. In this region,
\[
 \Ec(u(t))\ge\frac{1}{2}\Ac(u(t))-
 \frac{\Ac(u(t))^3}{6G^2}\ge\frac{1}{3}\Ac(u(t)),
\]
which proves \eqref{eq:trapping}. The $H^1$ bound follows from $c>0$.
\end{proof}

\section{Periodic profile analysis}\label{sec:periodic}

We collect the periodic estimates and profile results used below.
The linear and multilinear estimates are unchanged by the sign of
the quintic term; the hypotheses needed for the nonlinear Euclidean
input are verified in Section~\ref{sec:proof}.

\subsection{Local theory and stability}

We use the orthonormal Fourier basis
$(2\pi)^{-\frac{3}{2}}\e^{\ii x\cdot\xi}$, $\xi\in\Z^3$.
The atomic space $X^1$ is the restriction space associated with
\[
 \|u\|_{X^1(\R)}^2=
 \sum_{\xi\in\Z^3}\langle\xi\rangle^2
 \|\e^{\ii t|\xi|^2}\widehat u(t,\xi)\|_{U^2_t}^2.
\]
We use the standard definition of $U^2$ from \cite{HTT}. On intervals of
length at most one, this agrees with the convention in \cite{IP}; longer
intervals will be handled by finite gluing. In particular,
$X^1(I)$ controls $L_t^\infty H_x^1(I)$, and
$\|\e^{\ii t\Delta}f\|_{X^1(I)}\lesssim\|f\|_{H^1}$.
For any finite partition of $I$ into consecutive intervals $I_\nu$,
time truncation in $U^2$ and the triangle inequality give
\begin{equation}\label{eq:gluing}
 \|u\|_{X^1(I)}\lesssim\sum_{\nu=1}^{m}\|u\|_{X^1(I_\nu)}.
\end{equation}
For $p_0=\frac{41}{10}$ and $p_1=100$, set
\begin{equation}\label{eq:z}
 \|u\|_{Z(I)}=
 \sum_{p\in\{p_0,p_1\}}
 \sup_{\substack{J\subset I,  |J|\le1}}
 \Big(\sum_{N\ge1}N^{5-\frac{p}{2}}
 \|P_Nu\|_{L^p_{t,x}(J\times\T^3)}^p\Big)^{\frac{1}{p}}.
\end{equation}
We also write $\|u\|_{Z'(I)}=\|u\|_{Z(I)}^{\frac{1}{2}}\|u\|_{X^1(I)}^{\frac{1}{2}}$.
We denote the corresponding Duhamel forcing norm by $N(I)$; on an
interval containing $0$ we may take
\[
 \|e\|_{N(I)}=
 \Big\|\int_0^t\e^{\ii(t-s)\Delta}e(s)\,\mathrm ds\Big\|_{X^1(I)}.
\]
Changing the base time gives an equivalent norm.

We repeatedly use the quintilinear estimate \cite[Lemma 3.2]{IP}:
\begin{equation}\label{eq:quintilinear}
 \Big\|\prod_{k=1}^5\widetilde u_k\Big\|_{N(I)}
 \lesssim\sum_{k=1}^5\|u_k\|_{X^1(I)}
       \prod_{\ell\ne k}\|u_\ell\|_{Z'(I)},
 \qquad |I|\le1,
\end{equation}
where $\widetilde u_k$ is $u_k$ or $\overline{u_k}$.
We also use its frequency-ordered form \cite[equation (3.4)]{IP}:
for $D\ge1$ and $|I|\le1$,
\begin{equation}\label{eq:ordered}
 \Big\|\sum_{M\ge1}P_M\widetilde u_1
          \prod_{k=2}^{5}P_{\le DM}\widetilde u_k\Big\|_{N(I)}
 \lesssim_D \|u_1\|_{X^1(I)}
          \prod_{k=2}^{5}\|u_k\|_{Z'(I)}.
\end{equation}
Here $M$ is dyadic. This form allows us to specify which factor
carries the $X^1$ norm.

The proof of local well-posedness in \cite[Proposition 3.3]{IP}
uses these estimates and is unchanged for the focusing sign.
In the continuation argument, the $H^1$ bound is supplied by
Lemma~\ref{lem:trapping}.
We use the following consequences on intervals of length at most one:
\begin{enumerate}[label=(\roman*)]
 \item $\|u\|_{Z(I)}\lesssim\|u\|_{X^1(I)}$; a solution with bounded
 $H^1$ and finite $Z$ norm continues across a finite endpoint.
 \item A bound on $\|u\|_{L_t^\infty H_x^1(I)}+\|u\|_{Z(I)}$
 gives a bound on $\|u\|_{X^1(I)}$.
 \item There is $\varepsilon_*>0$ such that $\|f\|_{H^1}\le\varepsilon_*$
 gives a solution on $[-\frac{1}{2},\frac{1}{2}]$ satisfying
 \begin{equation}\label{eq:small-data}
  \|u\|_{X^1([-\frac{1}{2},\frac{1}{2}])}\lesssim\|f\|_{H^1},\qquad
  \|u-\e^{\ii t\Delta}f\|_{X^1([-\frac{1}{2},\frac{1}{2}])}
  \lesssim\|f\|_{H^1}^{5}.
 \end{equation}
\end{enumerate}
The last estimate follows from the quintic multilinear bound and the
small-data contraction argument. It is used only on a fixed finite
interval.

For later use, we state the stability conclusion of
\cite[Proposition 3.4]{IP}, which already includes both signs.
Suppose $|I|\le1$, $t_0\in I$, and
$\widetilde u\in C(I;H^1)\cap X^1(I)$ satisfies, in the Duhamel sense,
\[
 (\ii\partial_t+\Delta)\widetilde u=\F(\widetilde u)+e,
 \qquad
 \|\widetilde u\|_{L_t^\infty H_x^1(I)}
       +\|\widetilde u\|_{Z(I)}\le K,\qquad K\ge1.
\]
There are $\varepsilon_1(K)>0$ and $C(K)<\infty$ such that, if
\[
 \eta:=\|f-\widetilde u(t_0)\|_{H^1}+\|e\|_{N(I)}
       \le\varepsilon_1(K),
\]
then the solution with $u(t_0)=f$ exists on $I$ and satisfies
\begin{equation}\label{eq:stability}
 \|u\|_{X^1(I)}+\|\widetilde u\|_{X^1(I)}\le C(K),
 \qquad
 \|u-\widetilde u\|_{X^1(I)}\le C(K)\eta.
\end{equation}

\subsection{Linear profiles}

Fix a smooth cutoff $\chi$ supported in a coordinate ball and equal to
one near the origin. For large $N$, let
\begin{equation}\label{eq:embedding}
 T_N\phi(x)=N^{\frac{1}{2}}\chi(N^{\frac{1}{2}}x)\phi(Nx),
 \qquad \phi\in\dot H^1(\R^3),
\end{equation}
using the flat coordinate chart near the origin. The harmless fixed
chart size can be absorbed into $\chi$. These maps are uniformly
bounded from $\dot H^1$ to $H^1$, and
\begin{equation}\label{eq:embedding-limits}
 \|T_N\phi\|_2\to0,
 \quad \|\nabla T_N\phi\|_2^2\to \|\nabla\phi\|_2^2,
 \quad \|T_N\phi\|_6^6 \to \|\phi\|_6^6.
\end{equation}
These limits follow first for smooth compactly supported functions and
then by the uniform bounds and density in $\dot H^1$.

A Euclidean frame is a sequence $(N_n,t_n,x_n)$ with $N_n\to\infty$
and $t_n\to0$. Its linear profile is
\[
 p_n=\pi_{x_n}\e^{-\ii t_n\Delta}T_{N_n}\phi,
 \qquad \pi_y f(x)=f(x-y).
\]
Two frames are orthogonal if
\[
 \left|\log\frac{N_n^j}{N_n^k}\right|
 +(N_n^j)^2|t_n^j-t_n^k|
 +N_n^j\operatorname{dist}_{\T^3}(x_n^j,x_n^k)
 \to \infty.
\]
By \cite[Lemma 5.2]{IP}, equivalent frames may be normalized so that
either $t_n=0$ for every $n$, or $N_n^2t_n\to+\infty$ or $-\infty$.
We also normalize the frames jointly, as at the beginning of Case~III
in \cite[Section 6]{IP}. For each pair, the scales are either equal or
satisfy $\frac{N_n^j}{N_n^k}+\frac{N_n^k}{N_n^j}\to\infty$.
At equal scales, the time shifts are either equal or satisfy
$(N_n^j)^2|t_n^j-t_n^k|\to\infty$.
These equivalent replacements preserve the individual normalization;
the corresponding profiles are replaced according to
\cite[Lemma 5.2]{IP}. For each fixed $J$, the resulting $o_{H^1}(1)$
changes are absorbed into $r_n^J$ below. Finite free time shifts preserve
the gradient norm, but generally change the potential energy. All energy
calculations below therefore use the final normalized profiles and the
potential limits in Lemma~\ref{lem:time-escape}.

We use the linear decomposition of \cite[Proposition 5.4]{IP}.
If $f_n$ is bounded in $H^1$ and $I_n\ni0$ are intervals with
$|I_n|\to0$, then, after taking a subsequence,
\begin{equation}\label{eq:decomposition}
 f_n=g+\sum_{j=1}^{J}p_n^j+r_n^J,
\end{equation}
where $g\in H^1$, the Euclidean frames are pairwise orthogonal, and
\begin{gather}
 \Ac(f_n)=\Ac(g)+\sum_{j=1}^{J}\|\nabla\phi^j\|_2^2
                    +\Ac(r_n^J)+o_n(1),\label{eq:quadratic-split}\\
 \lim_{J\to\infty}\limsup_{n\to\infty}
 \Big|\|f_n\|_6^6-\|g\|_6^6-
                  \sum_{j=1}^{J}\|p_n^j\|_6^6\Big|=0,
                    \label{eq:potential-split}\\
 \lim_{J\to\infty}\limsup_{n\to\infty}
 \|\e^{\ii t\Delta}r_n^J\|_{Z(I_n)}=0.
                    \label{eq:small-remainder}
\end{gather}
The last assertion is the interpolation consequence used in
\cite[equation (6.3)]{IP}. Orthogonality also gives uniform $H^1$
bounds on finite sums of distinct profiles. We use
\eqref{eq:potential-split} in its stated order of limits.
For general intervals containing zero, we apply the cited proposition
to their symmetric hulls and then restrict to the original intervals.

\begin{lemma}\label{lem:time-escape}
For a normalized frame and its profile $p_n$,
\begin{equation}\label{eq:profile-potential}
 \lim_{n\to\infty}\|p_n\|_6^6=
 \begin{cases}
  \|\phi\|_{L^6(\R^3)}^6,&t_n=0,\\
  0,&N_n^2|t_n|\to \infty.
 \end{cases}
\end{equation}
\end{lemma}

\begin{proof}
The first assertion follows from \eqref{eq:embedding-limits}.
For the second, first let $\phi\in C_c^\infty(\R^3)$.
By \cite[equation (4.17)]{IP}, with $p=6$, for every fixed $T\ge1$,
\[
 \limsup_{n\to\infty}
 \|\e^{-\ii t_n\Delta}T_{N_n}\phi\|_6
 \lesssim_\phi T^{-\frac{1}{10}}.
\]
Indeed, $N_n^2|t_n|\to\infty$ and $t_n\to0$ imply
$TN_n^{-2}\le |t_n|\le T^{-1}$ for all sufficiently large $n$.
Let $T\to\infty$. The general case follows by density in $\dot H^1$,
the uniform boundedness of $T_N:\dot H^1\to H^1$, and the periodic
Sobolev embedding.
\end{proof}

\subsection{Nonlinear approximation and superposition}

We state the Euclidean approximation and superposition lemmas
in the form needed below, following \cite{IP}.

\begin{lemma}[Euclidean approximation]\label{lem:euclidean-approximation}
Let $p_n$ be associated with a normalized Euclidean frame and a profile
$\phi$. Suppose $V$ is a global focusing Euclidean solution with finite
$L^{10}_{t,x}$ norm, with $V(0)=\phi$ if $t_n=0$, and scattering to
$\phi$ as $t\to\mp\infty$ if $N_n^2t_n\to\pm\infty$.
For all sufficiently large $n$, the periodic solution $U_n$ with
$U_n(0)=p_n$ exists on a common interval $(-\tau,\tau)$, and
\begin{equation}\label{eq:nonlinear-profile}
 \limsup_{n\to\infty}\|U_n\|_{X^1(-\tau,\tau)}\le C_V<\infty.
\end{equation}
For every $0<\theta<1$, there is $T_0=T_0(V,\theta)$ such that
for each $T\ge T_0$ one can choose $R=R(V,\theta,T)\ge2T$
and obtain the following decomposition on
$L_T=(-T^{-1},T^{-1})\subset(-\tau,\tau)$, for all sufficiently large $n$.
Suppressing the dependence on $T$ and $R$, write
\begin{equation}\label{eq:core-decomposition}
 U_n=\omega_{n,\theta}+\omega_{n,\theta}^{-}
               +\omega_{n,\theta}^{+}+\rho_{n,\theta},
\end{equation}
where
\begin{align}
 \|\rho_{n,\theta}\|_{X^1(L_T)}
 +\sum_\pm\|\omega_{n,\theta}^{\pm}\|_{Z'(L_T)}
 &\le\theta,\label{eq:core-small}\\
 \|\omega_{n,\theta}\|_{X^1(L_T)}
 +\sum_\pm\|\omega_{n,\theta}^{\pm}\|_{X^1(L_T)}
 &\le C_V.\label{eq:core-bound}
\end{align}
The constant $C_V$ is independent of $\theta$ and $T$.
The core is supported where
\begin{equation}\label{eq:core-support}
 |t-t_n|\le TN_n^{-2},
 \qquad \operatorname{dist}_{\T^3}(x,x_n)\le RN_n^{-1}.
\end{equation}
In the interior of this time window, it satisfies
\begin{equation}\label{eq:core-derivatives}
 |\partial_x^\alpha\omega_{n,\theta}|
 +N_n^{-2}|\partial_t\partial_x^\alpha\omega_{n,\theta}|
 \lesssim_{V,\theta,T}N_n^{\frac{1}{2}+|\alpha|},
 \qquad |\alpha|\le6.
\end{equation}
For each fixed choice of these parameters,
\begin{equation}\label{eq:core-low}
 \lim_{B\to\infty}\limsup_{n\to\infty}
 \|P_{\le N_n/B}\omega_{n,\theta}\|_{X^1(L_T)}=0.
\end{equation}
Set
\[
 I_{n,T}^-=(-\infty,t_n-TN_n^{-2}),\qquad
 I_{n,T}^+=[t_n+TN_n^{-2},\infty).
\]
The tails have the form
\begin{equation}\label{eq:tail-structure}
 \omega_{n,\theta}^{\pm}(t)=
 \mathbf1_{I_{n,T}^{\pm}}(t)
 \e^{\ii(t-t_n)\Delta}\pi_{x_n}T_{N_n}\phi_\theta^\pm,
 \qquad t\in L_T.
\end{equation}
All time truncations use right-continuous representatives.
Here $\phi_\theta^\pm$ are regularized Euclidean scattering states satisfying
\[
 \|\phi_\theta^\pm\|_{\dot H^1(\R^3)}\le C_V,\qquad
 \|\phi_\theta^\pm\|_{L^1(\R^3)}\le R,\qquad
 \operatorname{supp}\widehat{\phi_\theta^\pm}
       \subset\{|\xi|\le R\}.
\]
In particular, a fixed finite collection of profiles can use a common
$T$ and, by enlargement, a common $R$.
\end{lemma}

\begin{proof}
We apply the Euclidean approximation argument of
\cite[Lemmas 4.2--4.3 and Proposition 4.4]{IP} to the
focusing Euclidean solution $V$. Its finite $L^{10}_{t,x}$
norm and the Euclidean local theory give
$\nabla V\in L^{10/3}_{t,x}(\R\times\R^3)$. The truncation
and stability arguments are unchanged for the focusing sign,
while the extinction estimate is linear; cf.
\cite[Proposition 5.8]{Yue}. Together with scattering, this
gives $\tau=\tau(V)>0$ such that, for all sufficiently large
$n$, the solution $U_n$ exists on $(-\tau,\tau)$ and satisfies
\eqref{eq:nonlinear-profile}.

Fix $0<\theta<1$. The construction leading to
\cite[equation (7.8)]{IP}, cf. \cite[Corollary 5.9]{Yue},
gives $T_0=T_0(V,\theta)$ such that, for every $T\ge T_0$,
one may choose $R=R(V,\theta,T)\ge2T$ so that
\eqref{eq:core-decomposition}, \eqref{eq:core-small},
\eqref{eq:core-support}--\eqref{eq:core-derivatives}, and
\eqref{eq:tail-structure} hold. After regularizing the
scattering states, absorbing the resulting errors into
$\rho_{n,\theta}$, and enlarging $R$ if necessary, we may
also arrange that
\begin{equation*}
 \|\phi_\theta^\pm\|_{\dot H^1(\R^3)}\le C_V,\qquad
 \|\phi_\theta^\pm\|_{L^1(\R^3)}\le R,\qquad
 \operatorname{supp}\widehat{\phi_\theta^\pm}
 \subset\{|\xi|\le R\}.
\end{equation*}

We prove \eqref{eq:core-low}. Fix $\theta,T$ and the
corresponding $R$, and set
\begin{equation*}
 N=N_n,\qquad \omega=\omega_{n,\theta},\qquad
 J_n=L_T\cap[t_n-TN^{-2},t_n+TN^{-2}).
\end{equation*}
By \eqref{eq:core-support}--\eqref{eq:core-derivatives} and
Bernstein's inequality, for $1\le K\le N$, we have
\begin{align*}
 \sup_{t\in J_n}\|P_{\le K}\omega(t)\|_{H^1_x}
 &\lesssim_{V,\theta,T}K^{5/2}N^{-5/2},\\
 \sup_{t\in J_n}
 \bigl(
  \|P_{\le K}\partial_t\omega(t)\|_{H^1_x}
  +\|\Delta P_{\le K}\omega(t)\|_{H^1_x}
 \bigr)
 &\lesssim_{V,\theta,T}
 K^{5/2}N^{-1/2}+K^{9/2}N^{-5/2}.
\end{align*}
The time-translated form of \cite[equation (2.5)]{IP},
together with the boundedness of right-continuous interval
truncations in $U^2_\Delta$, therefore gives
\begin{align*}
\|P_{\le K}\omega\|_{X^1(L_T)}
&\lesssim \sup_{t\in J_n}\|P_{\le K}\omega(t)\|_{H^1_x}
+\int_{J_n}\|(i\partial_t+\Delta)P_{\le K}\omega(t)\|_{H^1_x}\,\mathrm{d}t\\
&\lesssim_{V,\theta,T}
K^{5/2}N^{-5/2}+K^{9/2}N^{-9/2}.
\end{align*}
where $|J_n|\le2TN^{-2}$ was used. For fixed $B\ge1$, take
$K=N_n/B$ according to the convention for arbitrary frequency
thresholds. Since $1\le K\le N_n$ for all sufficiently large
$n$, it follows that
\begin{equation*}
 \limsup_{n\to\infty}
 \|P_{\le N_n/B}\omega_{n,\theta}\|_{X^1(L_T)}
 \lesssim_{V,\theta,T}
 B^{-5/2}+B^{-9/2}.
\end{equation*}
Letting $B\to\infty$ proves \eqref{eq:core-low}.

Finally, \eqref{eq:tail-structure}, the interval-truncation
property of $X^1$, and the scaling limits for $T_N$ give
\begin{equation*}
 \limsup_{n\to\infty}
 \sum_\pm\|\omega_{n,\theta}^\pm\|_{X^1(L_T)}
 \lesssim
 \sum_\pm\|\phi_\theta^\pm\|_{\dot H^1(\R^3)}
 \le C_V.
\end{equation*}
Moreover, \eqref{eq:core-decomposition} implies
\begin{align*}
 &\|\omega_{n,\theta}\|_{X^1(L_T)}
 +\sum_\pm\|\omega_{n,\theta}^\pm\|_{X^1(L_T)}\\
 &\qquad\le
 \|U_n\|_{X^1(L_T)}
 +2\sum_\pm\|\omega_{n,\theta}^\pm\|_{X^1(L_T)}
 +\|\rho_{n,\theta}\|_{X^1(L_T)}.
\end{align*}
Thus \eqref{eq:nonlinear-profile} and
\eqref{eq:core-small} yield \eqref{eq:core-bound}, after
enlarging $C_V$ by a fixed factor. The resulting constant
depends only on $V$ and is independent of $\theta$ and $T$.

For a fixed finite collection of profiles, take $T$ above
all the corresponding thresholds and then choose a common
sufficiently large $R$.
\end{proof}

We conclude this section with the superposition estimate needed
to control the nonlinear error of the approximate solution.

\begin{lemma}\label{lem:superposition}
Let $L$ be a fixed interval containing zero, with $|L|\le1$, and let
$I_n\subset L$, $0\in I_n$, satisfy $|I_n|\to0$.
Let $v\in X^1(L)$ be a fixed periodic solution, and let $U_n^j$
be nonlinear profiles as in Lemma~\ref{lem:euclidean-approximation},
associated with pairwise orthogonal frames normalized jointly as above,
and defined on $L$ for all
sufficiently large $n$ depending on $j$.
Set
\[
 S_n^J=v+\sum_{j=1}^{J}U_n^j,\qquad
 R_n^J=\e^{\ii t\Delta}r_n^J.
\]
Suppose
\begin{gather}
 \sup_J\limsup_{n\to\infty}
 \bigl(\|S_n^J\|_{X^1(L)}+\|R_n^J\|_{X^1(L)}\bigr)\le Q<\infty,
 \label{eq:assembly-hypothesis}\\
 \lim_{A\to\infty}\sup_{J>A}\limsup_{n\to\infty}
 \|S_n^J-S_n^A\|_{X^1(L)}=0,\label{eq:cauchy-hypothesis}\\
 \lim_{J\to\infty}\limsup_{n\to\infty}
 \|R_n^J\|_{Z(I_n)}=0.\label{eq:remainder-hypothesis}
\end{gather}
Then the first limit below holds for each fixed $J$, and the second
holds in the indicated order:
\begin{gather}
 \limsup_{n\to\infty}
 \|\F(S_n^J)-\F(v)-\sum_{j=1}^{J}\F(U_n^j)\|_{N(I_n)}=0,
 \label{eq:interaction-one}\\
 \lim_{J\to\infty}\limsup_{n\to\infty}
 \|\F(S_n^J+R_n^J)-\F(S_n^J)\|_{N(I_n)}=0.
 \label{eq:interaction-two}
\end{gather}
\end{lemma}

\begin{proof}
We follow the proof of \cite[Lemma 6.2]{IP}; cf.
\cite[Lemma 7.4 and Section 8]{Yue}. As in \cite{IP},
$\mathfrak O_{4,1}(a,b)$ denotes a quintilinear expression
containing four factors from $\{a,\overline a\}$ and one
factor from $\{b,\overline b\}$. 
All other products below may likewise contain complex conjugates.
By \eqref{eq:assembly-hypothesis},
\eqref{eq:remainder-hypothesis}, and the definition of $Z'$,
\begin{equation*}
 \lim_{J\to\infty}\limsup_{n\to\infty}
 \|R_n^J\|_{Z'(I_n)}=0.
\end{equation*}
Since $v$ is fixed and $|I_n|\to0$, we also have
\begin{equation*}
 \|v\|_{Z'(I_n)}\to 0,
 \qquad
 \|P_{>K}v\|_{X^1(L)}\to 0
 \quad\text{as }K\to\infty.
\end{equation*}

We first prove \eqref{eq:interaction-one}. Fix $J$ and
apply \eqref{eq:core-decomposition} to the first $J$ profiles,
with common parameters $T,R$. For fixed $\theta$, we have
$I_n\subset L_T$ for all sufficiently large $n$.
Expanding the quintic interaction, terms containing an error
or at least two background or scattering factors are
controlled by \eqref{eq:quintilinear}. Terms containing
two distinct cores tend to zero by \cite[Lemma 7.2]{IP}.
The remaining core--tail terms satisfy
\begin{equation*}
 \limsup_{n\to\infty}
 \|\mathfrak O_{4,1}
   (\omega_{n,\theta}^j,\omega_{n,\theta}^{k,\pm})
 \|_{N(I_n)}
 \lesssim_J\theta^{1/10},
 \qquad j\ne k,
\end{equation*}
by the argument proving \cite[equation (7.11)]{IP}.
Its four cases apply by the joint normalization:
separated scales are treated by
\cite[Lemma 7.1]{IP} and \eqref{eq:ordered};
equal scales and equal time centers give disjoint time supports;
the remaining case follows by smooth approximation of the
scattering state and \cite[equation (4.17)]{IP}.

For $\mathfrak O_{4,1}(\omega_{n,\theta}^j,v)$, first
discard $P_{>K}v$ and
$P_{\le N_n^j/B}\omega_{n,\theta}^j$ using
\eqref{eq:quintilinear}. For fixed $K,B$, the remaining
core frequencies exceed those of $P_{\le K}v$ when $n$
is sufficiently large. Thus \eqref{eq:ordered} places
the $X^1$ norm on a core factor and gives a factor
$\|v\|_{Z'(I_n)}\to0$. Letting $B\to\infty$ and then
$K\to\infty$, using \eqref{eq:core-low}, proves that
these terms also vanish. Hence
\begin{equation*}
 \limsup_{n\to\infty}
 \Big\|
 \F(S_n^J)-\F(v)-\sum_{j=1}^J\F(U_n^j)
 \Big\|_{N(I_n)}
 \lesssim_J\theta^{1/10}.
\end{equation*}
Letting $\theta\to0$ proves \eqref{eq:interaction-one}.

We next prove \eqref{eq:interaction-two}. Set $ \delta_A
 :=
 \sup_{J>A}\limsup_{n\to\infty}
 \|S_n^J-S_n^A\|_{X^1(L)}.$ 
By \eqref{eq:cauchy-hypothesis}, $\delta_A\to0$.
Expanding the quintic difference and applying
\eqref{eq:quintilinear}, we obtain, for $J>A$,
\begin{align*}
 &\|\F(S_n^J+R_n^J)-\F(S_n^J)\|_{N(I_n)}\\
 &\quad\lesssim_Q
 \|R_n^J\|_{Z'(I_n)}
 +\|S_n^J-S_n^A\|_{X^1(L)}
 +\sum_{\mathfrak O_{4,1}}
  \|\mathfrak O_{4,1}(S_n^A,R_n^J)\|_{N(I_n)},
\end{align*}
where the sum runs over the finitely many conjugation patterns.

Fix $A$ and decompose its profiles with common parameters
$\theta,T,R$. Terms containing an error contribute
$O_A(\theta)$. Terms containing a background or scattering
factor contribute \( O_A\bigl(
 \theta+\|v\|_{Z'(I_n)}+\|R_n^J\|_{Z'(I_n)}
 \bigr)\)
by \eqref{eq:quintilinear}. Terms containing two distinct
cores vanish as $n\to\infty$ by \cite[Lemma 7.2]{IP}.
It remains to estimate
$\mathfrak O_{4,1}(\omega_{n,\theta}^j,R_n^J)$.
Write $N=N_n^j$ and $\omega=\omega_{n,\theta}^j$.
The localized form of \cite[Lemma 7.1]{IP}, together with
\eqref{eq:assembly-hypothesis}, gives
\begin{equation*}
 \limsup_{n\to\infty}
 \|\mathfrak O_{4,1}(\omega,P_{>BN}R_n^J)\|_{N(I_n)}
 \lesssim_{A,\theta,T,R,Q}B^{-1/200},
\end{equation*}
with a constant independent of $J$. Here
$\limsup_n\|r_n^J\|_{H^1}\lesssim Q$, since
$R_n^J(0)=r_n^J$.

Define
\begin{equation*}
 \eta_{A,\theta}(B)
 :=
 \max_{1\le j\le A}\limsup_{n\to\infty}
 \|P_{\le N_n^j/B}\omega_{n,\theta}^j\|_{X^1(I_n)}.
\end{equation*}
By \eqref{eq:core-low}, $\eta_{A,\theta}(B)\to0$ as
$B\to\infty$. Discarding these low-frequency core factors
costs $O_{A,Q}(\eta_{A,\theta}(B))$ by
\eqref{eq:quintilinear}. For the remaining product,
the largest core frequency $M$ satisfies $M\gtrsim N/B$,
whereas the remainder frequency is $O(BN)$.
Thus \eqref{eq:ordered}, with $D\lesssim B^2$, yields
\begin{equation*}
 \|\mathfrak O_{4,1}
   (P_{>N/B}\omega,P_{\le BN}R_n^J)\|_{N(I_n)}
 \lesssim_B
 \|\omega\|_{X^1(I_n)}
 \|\omega\|_{Z'(I_n)}^3
 \|R_n^J\|_{Z'(I_n)}.
\end{equation*}
The core bounds are independent of $J$. Consequently,
\begin{align*}
 &\limsup_{J\to\infty}\limsup_{n\to\infty}
 \|\F(S_n^J+R_n^J)-\F(S_n^J)\|_{N(I_n)}\\
 &\quad\lesssim_Q\delta_A
 +C_{A,Q}\theta
 +C_{A,\theta,T,R,Q}B^{-1/200}
 +C_{A,Q}\eta_{A,\theta}(B).
\end{align*}
Letting $B\to\infty$, then $\theta\to0$, and finally
$A\to\infty$ proves \eqref{eq:interaction-two}.
\end{proof}

\section{Proof of Theorem~\ref{thm:main}}\label{sec:proof}

\subsection{Euclidean profile energies}

\begin{lemma}\label{lem:profile-energy}
Suppose the sequence in \eqref{eq:decomposition} satisfies
\[
 \Ac(f_n)\le A_*<G,
 \qquad \Ec(f_n)\le E_*<\frac{G}{3}.
\]
For each profile define
\begin{equation}\label{eq:profile-energy}
 e_j=
 \begin{cases}
  \ER(\phi^j),&t_n^j=0,\\
  \frac{1}{2}\|\nabla\phi^j\|_2^2,
      &(N_n^j)^2|t_n^j|\to \infty.
 \end{cases}
\end{equation}
Then
\begin{equation}\label{eq:profile-threshold}
 \|\nabla\phi^j\|_2^2\le A_*<G,
 \qquad 0\le e_j\le E_*<\frac{G}{3}.
\end{equation}
Each profile therefore determines a global Euclidean focusing
solution with finite $L^{10}_{t,x}$ norm and the matching condition
needed in Lemma~\ref{lem:euclidean-approximation}.
\end{lemma}

\begin{proof}
Passing to a subsequence, assume that
\begin{equation*}
 \Ac(f_n)\to q,\qquad
 \|f_n\|_6^6\to b,\qquad
 \Ec(f_n)\to e=\frac12q-\frac16b\le E_*.
\end{equation*}
Set $a_j=\|\nabla\phi^j\|_2^2$ and
$b_j=\lim_{n\to\infty}\|p_n^j\|_6^6$; the latter limit is
given by Lemma~\ref{lem:time-escape}. For each fixed $J$,
\eqref{eq:quadratic-split} gives
\begin{equation*}
 d_J:=\lim_{n\to\infty}\Ac(r_n^J)
 =q-\Ac(g)-\sum_{j=1}^J a_j\ge0.
\end{equation*}
Thus $d_J\downarrow d\ge0$. Taking $n\to\infty$ first
and then $J\to\infty$ in \eqref{eq:potential-split},
we obtain
\begin{equation}\label{eq:complete-budget}
 q=\Ac(g)+\sum_{j\ge1}a_j+d,
 \qquad
 b=\|g\|_6^6+\sum_{j\ge1}b_j.
\end{equation}
In particular, $\Ac(g)\le A_*<G$ and $a_j\le A_*<G$.
The sharp Sobolev inequalities imply
\begin{equation*}
 \|g\|_6^6\le G^{-2}\Ac(g)^3,
 \qquad
 0\le b_j\le G^{-2}a_j^3,
\end{equation*}
and hence
\begin{equation*}
 \Ec(g)\ge\frac13\Ac(g)\ge0,
 \qquad
 e_j=\frac12a_j-\frac16b_j\ge\frac13a_j\ge0.
\end{equation*}
Consequently, \eqref{eq:complete-budget} yields
\begin{equation}\label{eq:positive-budget}
 e=\Ec(g)+\sum_{j\ge1}e_j+\frac d2.
\end{equation}
This proves \eqref{eq:profile-threshold}.

If $t_n^j=0$, apply Theorem~\ref{thm:euclidean} to
$\phi^j$. Otherwise, the Euclidean wave operator
\cite[Remark 2.13]{KM} gives a solution $V^j$ on a time
tail scattering to $\phi^j$ in the prescribed direction.
By density and the dispersive estimate,
$\e^{\ii t\Delta_{\R^3}}\phi^j\to0$ in $L^6$ as
$t\to\pm\infty$. Scattering, Sobolev embedding, and energy
conservation therefore give
\begin{equation*}
 \ER(V^j)=e_j=\frac12a_j<\frac G3,
 \qquad
 \|\nabla V^j(t)\|_2^2\to  a_j
 =2e_j<\frac{2G}{3}<G
\end{equation*}
in the prescribed time direction. At a sufficiently distant
finite time, both hypotheses of
Theorem~\ref{thm:euclidean} hold. Thus $V^j$ extends to a
global scattering solution with the required matching condition.
\end{proof}

\subsection{Uniform bounds on short intervals}
We now establish uniform short-time spacetime bounds by contradiction. 

\begin{proposition}\label{prop:short}
For each $0<E_0<\frac{G}{3}$, there exist
$0<\tau_0\le1$ and $B_0<\infty$, depending only on
$E_0$ and $c$, such that every solution $u$ with
$\Ac(u_0)<G$ and $\Ec(u_0)\le E_0$ satisfies
\begin{equation}\label{eq:short}
 \|u\|_{Z(I)}\le B_0
\end{equation}
on every compact interval $I$ in its lifespan with
$|I|\le\tau_0$.
\end{proposition}

\begin{proof}
Suppose otherwise. After time translation, choose solutions
$u_n$ and compact intervals $I_n\ni0$ in their lifespans
such that
\begin{equation}\label{eq:bad-sequence}
 |I_n|\to 0,
 \qquad
 \|u_n\|_{Z(I_n)}\to \infty.
\end{equation}
By Lemma~\ref{lem:trapping},
$\Ac(u_n(t))\le3E_0<G$. Apply
\eqref{eq:decomposition}--\eqref{eq:small-remainder}
to $f_n=u_n(0)$. Lemma~\ref{lem:profile-energy}, with
$A_*=3E_0$ and $E_*=E_0$, supplies the Euclidean solutions
needed for the nonlinear profiles $U_n^j$.
Let $v$ be the local periodic solution with data $g$, and set \(S_n^J=v+\sum_{j=1}^J U_n^j\).

We first choose a common interval for the nonlinear profiles.
Write $a_j=\|\nabla\phi^j\|_2^2$, discarding zero profiles.
By \eqref{eq:quadratic-split}, $\sum_j a_j<\infty$.
The scaling limits for $T_N$ give
$\|p_n^j\|_{H^1}^2\to a_j$ for each fixed $j$.
Choose $J_0$ so that $2a_j^{1/2}$ is below the small-data
threshold in \eqref{eq:small-data} whenever $j>J_0$.
For the finitely many profiles $j\le J_0$, use
Lemma~\ref{lem:euclidean-approximation}.
Choose a fixed interval $L=[-\tau,\tau]$,
$0<\tau\le1/2$, contained in their common existence
intervals and in the local lifespan of $v$.
Eventually $I_n\subset L$.

For each fixed $J>J_0$ and all sufficiently large $n$,
\eqref{eq:small-data} gives, for $J_0<j\le J$,
\begin{equation}\label{eq:tail}
 \|U_n^j\|_{X^1(L)}\lesssim a_j^{1/2},
 \qquad
 \|U_n^j-\e^{\ii t\Delta}p_n^j\|_{X^1(L)}
 \lesssim a_j^{5/2}.
\end{equation}
For fixed $J>A\ge J_0$, write
\begin{equation*}
 S_n^J-S_n^A
 =
 \e^{\ii t\Delta}\sum_{A<j\le J}p_n^j
 +
 \sum_{A<j\le J}
 \bigl(U_n^j-\e^{\ii t\Delta}p_n^j\bigr).
\end{equation*}
The linear estimate, orthogonality, and \eqref{eq:tail}
therefore imply
\begin{equation}\label{eq:tail-cauchy}
 \sup_{J>A}\limsup_{n\to\infty}
 \|S_n^J-S_n^A\|_{X^1(L)}
 \lesssim
 \Big(\sum_{j>A}a_j\Big)^{1/2}
 +\sum_{j>A}a_j^{5/2}
 \to 0
\end{equation}
as $A\to\infty$. Combining the same bounds with
\eqref{eq:nonlinear-profile} for $j\le J_0$ and the
quadratic bound on $r_n^J$, we obtain some $Q\ge1$ such that
\begin{equation}\label{eq:uniform-assembly}
 \sup_J\limsup_{n\to\infty}
 \bigl(
 \|S_n^J\|_{X^1(L)}
 +\|\e^{\ii t\Delta}r_n^J\|_{X^1(L)}
 \bigr)\le Q.
\end{equation}
Here $Q$ may depend on the extracted profiles but is
independent of $J$. Thus
\eqref{eq:tail-cauchy}, \eqref{eq:uniform-assembly},
and \eqref{eq:small-remainder} verify the hypotheses of
Lemma~\ref{lem:superposition}.

Set
\begin{equation*}
 u_{n,\mathrm{app}}^J=S_n^J+\e^{\ii t\Delta}r_n^J.
\end{equation*}
By \eqref{eq:decomposition},
$u_{n,\mathrm{app}}^J(0)=u_n(0)$. Its equation error is
\begin{equation*}
 e_n^J:=
 (i\partial_t+\Delta)u_{n,\mathrm{app}}^J
 -\F(u_{n,\mathrm{app}}^J)
 =
 \F(v)+\sum_{j=1}^J\F(U_n^j)
 -\F(u_{n,\mathrm{app}}^J).
\end{equation*}
Lemma~\ref{lem:superposition} gives
\begin{equation}\label{eq:error}
 \lim_{J\to\infty}\limsup_{n\to\infty}
 \|e_n^J\|_{N(I_n)}=0.
\end{equation}
Moreover, \eqref{eq:uniform-assembly} and
$Z\lesssim X^1$ give a constant $K=CQ\ge1$ such that,
for each fixed $J$ and all sufficiently large $n$,
\begin{equation*}
 \|u_{n,\mathrm{app}}^J\|_{L_t^\infty H_x^1(I_n)}
 +\|u_{n,\mathrm{app}}^J\|_{Z(I_n)}
 \le K.
\end{equation*}
Choose $J$ so large that the limsup in
\eqref{eq:error} is less than
$\varepsilon_1(K)/2$, where $\varepsilon_1(K)$ is the
threshold in \eqref{eq:stability}. Fix this $J$ and take
$n$ sufficiently large. Stability applies with zero
initial discrepancy; by uniqueness, the resulting solution
agrees with $u_n$ on $I_n$. Hence
\begin{equation*}
 \|u_n\|_{Z(I_n)}
 \lesssim\|u_n\|_{X^1(I_n)}
 \le C(K),
\end{equation*}
contradicting \eqref{eq:bad-sequence}.
\end{proof}

\subsection{Global existence and continuous dependence}

\begin{proof}[Completion of the proof of Theorem~\ref{thm:main}]
Choose $0<E_0<\frac{G}{3}$ with $\Ec(u_0)<E_0$. Energy trapping bounds
the $H^1$ norm throughout the maximal lifespan. If a positive
endpoint $T_+$ were finite, Proposition~\ref{prop:short} would bound
the $Z$ norm uniformly on compact subintervals of
$(T_+-\frac{\tau_0}{2},T_+)$. Passing to the endpoint and using the
continuation criterion contradicts maximality. The negative
endpoint is treated by time reversal.

For the uniform assertion, partition $[-T,T]$ into a bounded number
of intervals of length at most $\tau_0$. On each interval,
Proposition~\ref{prop:short}, the $H^1$ bound, and the local
estimate in Section~\ref{sec:periodic} give a common $X^1$ bound.
Applying \eqref{eq:gluing} gives \eqref{eq:uniform-main}.

Finally, the strict inequalities defining the initial-data class
persist in a sufficiently small $H^1$ neighborhood of each datum.
Take the reference solution as the approximate solution in
\eqref{eq:stability}, with zero forcing error. Starting at time zero
and proceeding in both directions, stability propagates convergence
through the finitely many partition points.
The gluing estimate \eqref{eq:gluing} then gives convergence in
$X^1([-T,T])$. Uniqueness follows from the local critical theory.
\end{proof}

\begin{remark}
For the corresponding defocusing result on rectangular tori, see
Strunk~\cite{Strunk}. The present proof is restricted to the standard torus.
\end{remark}

\section*{Declarations}

\subsection*{Acknowledgments}
M. Chen and H. Wang were supported by the National Key R\&D Program of China
(2025YFA1018500) and the Beijing Natural Science Foundation (No.~1262019).

\subsection*{Data availability}
No datasets were generated or analyzed for this article.

\subsection*{Competing interests}
The authors declare no competing interests.

\subsection*{Use of AI}
ChatGPT (OpenAI) assisted with mathematical discussions and the
drafting and revision of parts of the proofs. It also assisted with
language editing, LaTeX preparation, bibliographic searches, and
consistency checks.
The authors assume responsibility for all content.


\begin{thebibliography}{99}

\bibitem{Aubin}
T.~Aubin,
\emph{Probl\`emes isop\'erim\'etriques et espaces de Sobolev},
J. Differential Geom. \textbf{11} (1976), no.~4, 573--598.
\href{https://doi.org/10.4310/jdg/1214433725}{doi:10.4310/jdg/1214433725}.

\bibitem{CH}
P.-h.~Chung and D.~Han,
\emph{Threshold scattering for the three-dimensional focusing
energy-critical nonlinear Schr\"odinger equation}, preprint (2026),
\href{https://arxiv.org/abs/2610.10658v1}{arXiv:2610.10658v1}.

\bibitem{CKSTT}
J.~Colliander, M.~Keel, G.~Staffilani, H.~Takaoka, and T.~Tao,
\emph{Global well-posedness and scattering for the energy-critical
nonlinear Schr\"odinger equation in $\R^3$},
Ann. of Math. (2) \textbf{167} (2008), no.~3, 767--865.
\href{https://doi.org/10.4007/annals.2008.167.767}{doi:10.4007/annals.2008.167.767}.

\bibitem{Dodson}
B.~Dodson,
\emph{Global well-posedness and scattering for the focusing, cubic
Schr\"odinger equation in dimension $d=4$},
Ann. Sci. \'{E}c. Norm. Sup\'{e}r. (4) \textbf{52} (2019), no.~1, 139--180.
\href{https://doi.org/10.24033/asens.2385}{doi:10.24033/asens.2385}.

\bibitem{HTT}
S.~Herr, D.~Tataru, and N.~Tzvetkov,
\emph{Global well-posedness of the energy-critical nonlinear
Schr\"odinger equation with small initial data in $H^1(\T^3)$},
Duke Math. J. \textbf{159} (2011), no.~2, 329--349.
\href{https://doi.org/10.1215/00127094-1415889}{doi:10.1215/00127094-1415889}.

\bibitem{IP}
A.~D. Ionescu and B.~Pausader,
\emph{The energy-critical defocusing NLS on $\T^3$},
Duke Math. J. \textbf{161} (2012), no.~8, 1581--1612.
\href{https://doi.org/10.1215/00127094-1593335}{doi:10.1215/00127094-1593335}.
Numbered statements are cited from
\href{https://arxiv.org/abs/1102.5771v1}{arXiv:1102.5771v1}.

\bibitem{KM}
C.~E. Kenig and F.~Merle,
\emph{Global well-posedness, scattering and blow-up for the
energy-critical, focusing, non-linear Schr\"odinger equation in the
radial case},
Invent. Math. \textbf{166} (2006), no.~3, 645--675.
\href{https://doi.org/10.1007/s00222-006-0011-4}{doi:10.1007/s00222-006-0011-4}.
Numbered statements are cited from
\href{https://arxiv.org/abs/math/0610266v1}{arXiv:math/0610266v1}.

\bibitem{KV}
R.~Killip and M.~Visan,
\emph{The focusing energy-critical nonlinear Schr\"odinger equation
in dimensions five and higher},
Amer. J. Math. \textbf{132} (2010), no.~2, 361--424.
\href{https://doi.org/10.1353/ajm.0.0107}{doi:10.1353/ajm.0.0107}.

\bibitem{Kwak}
B.~Kwak,
\emph{Global well-posedness of the energy-critical nonlinear
Schr\"odinger equations on $\T^d$},
Ann. PDE \textbf{12} (2026), article~19.
\href{https://doi.org/10.1007/s40818-026-00246-2}{doi:10.1007/s40818-026-00246-2}.

\bibitem{Luo}
Y.~Luo,
\emph{On long time behavior of the focusing energy-critical NLS on
$\R^d\times\T$ via semivirial-vanishing geometry},
J. Math. Pures Appl. (9) \textbf{177} (2023), 415--454.
\href{https://doi.org/10.1016/j.matpur.2023.07.006}{doi:10.1016/j.matpur.\allowbreak2023.07.006}.

\bibitem{Strunk}
N.~Strunk,
\emph{Global well-posedness of the energy-critical defocusing NLS
on rectangular tori in three dimensions},
Differential Integral Equations \textbf{28} (2015), no.~11--12,
1069--1084.
\href{https://arxiv.org/abs/1411.7886}{arXiv:1411.7886}.

\bibitem{SZ}
Q.~Su and Z.~Zhao,
\emph{Scattering below the ground state for the three-dimensional
focusing energy-critical NLS}, preprint (2026),
\href{https://arxiv.org/abs/2609.30228v3}{arXiv:2609.30228v3}.

\bibitem{Talenti}
G.~Talenti,
\emph{Best constant in Sobolev inequality},
Ann. Mat. Pura Appl. (4) \textbf{110} (1976), 353--372.
\href{https://doi.org/10.1007/BF02418013}{doi:10.1007/BF02418013}.

\bibitem{YYZ}
X.~Yu, H.~Yue, and Z.~Zhao,
\emph{Global well-posedness for the focusing cubic NLS on the
product space $\R\times\T^3$},
SIAM J. Math. Anal. \textbf{53} (2021), no.~2, 2243--2274.
\href{https://doi.org/10.1137/20M1364953}{doi:10.1137/20M1364953}.

\bibitem{Yue}
H.~Yue,
\emph{Global well-posedness for the energy-critical focusing
nonlinear Schr\"odinger equation on $\T^4$},
J. Differential Equations \textbf{280} (2021), 754--804.
\href{https://doi.org/10.1016/j.jde.2021.01.031}{doi:10.1016/j.jde.2021.01.031}.
The defocusing result is stated as Theorem~1.1 in the earlier version
\href{https://arxiv.org/abs/1805.09816v1}{arXiv:1805.09816v1}.

\end{thebibliography}
\end{document}